\documentclass[11pt,reqno]{amsart}
\usepackage[T1]{fontenc}
\usepackage{lmodern,microtype}
\usepackage[a4paper,margin=28mm]{geometry}
\usepackage{amsmath,amssymb,amsthm,mathtools}
\usepackage{needspace}
\usepackage[colorlinks=true,linkcolor=blue,citecolor=blue,urlcolor=blue]{hyperref}
\newtheorem{maintheorem}{Theorem}

\newtheorem{maincorollary}[maintheorem]{Corollary}
\newtheorem{theorem}{Theorem}[section]
\newtheorem{lemma}[theorem]{Lemma}
\newtheorem{proposition}[theorem]{Proposition}

\theoremstyle{definition}
\newtheorem{definition}[theorem]{Definition}
\theoremstyle{remark}
\newtheorem{remark}[theorem]{Remark}
\numberwithin{equation}{section}
\newcommand{\Cr}{C_r^*}
\newcommand{\PHP}{\mathrm P_{\mathrm{PHP}}}
\newcommand{\ind}{\mathbf 1}
\newcommand{\can}{\operatorname{can}}
\newcommand{\red}{\operatorname{red}}

\title[$C^*$-selflessness for free Burnside groups]{$C^*$-selflessness for free Burnside groups}
\date{}
\author{Anna Cascioli}
\email[Anna Cascioli]{acasciol@uni-muenster.de}
\address{University of Münster, Einsteinstrasse 62, Münster 48149, Germany}

\begin{document}
\begin{abstract}
The free Burnside group $B(m,n)$ is obtained from the free group of rank $m$ by imposing the relations $w^n=1$ for all words $w$. We prove that $B(m,n)$ is completely $C^*$-selfless for every $m\geq2$ and every odd exponent $n\geq557$. Consequently, its reduced group $C^*$-algebra has stable rank one and strict comparison of positive elements by its
canonical trace. This answers the question of whether $C^*$-selfless groups can satisfy a nontrivial group law.
The proof establishes Ozawa's property $\PHP$ by combining a combinatorial construction with a relative power-avoidance
argument based on the canonical forms of Atkarskaya--Rips--Tent.
\end{abstract}
\maketitle

\section{Introduction}

Since its introduction by Robert \cite{Robert}, selflessness has attracted substantial interest as a source of regularity for $C^*$-algebras beyond the nuclear setting. It expresses the possibility of adjoining a freely independent copy of an algebra inside a norm ultrapower while keeping the diagonal copy fixed.
In the presence of a faithful tracial state, selflessness implies strict comparison and other structural properties.
Nonabelian free groups provide basic examples in the setting of reduced group $C^*$-algebras \cite{Robert,AGKEP}. Ozawa \cite{Ozawa} developed dynamical and combinatorial criteria for the stronger property of \emph{complete selflessness}.

We study whether this free-product phenomenon can coexist with a law of bounded exponent. Let $F_m$ be the free group of rank $m\geq2$, allowing $m$ to be infinite. For an integer $n\geq2$, the free Burnside group $B(m,n)$ is defined by
\[
 B(m,n)=F_m/\langle\!\langle w^n:w\in F_m\rangle\!\rangle.
\] Every element has order dividing $n$, so these groups are torsion and contain no nonabelian free subgroup.
For finite $m\geq2$, Novikov and Adian \cite{NovikovAdian1968I, NovikovAdian1968II, NovikovAdian1968III} established infinitude for sufficiently large odd exponents. Adian \cite{Adian} subsequently obtained the bound $n\geq665$ and later proved nonamenability for the same odd exponents \cite{Adian1983}.
Atkarskaya, Rips, and Tent \cite{ART} lowered the infinitude
threshold to $557$, the best currently known uniform threshold
for odd exponents. The present paper builds on the techniques developed in \cite{Nonamenability} for a direct proof of nonamenability in this exponent range.

Olshanskii and Osin \cite[Theorem~1.2]{OO} proved simplicity and uniqueness of trace for $\Cr(B(m,n))$ in rank $m\geq2$ and sufficiently large odd exponent. Patchell
\cite[p.~1909]{Patchell} later raised the question of selflessness for free Burnside groups. We answer this question throughout the exponent range of \cite{ART}, uniformly in the rank. More precisely, we verify Ozawa's combinatorial property $\PHP$, which implies complete selflessness \cite[Theorem~14]{Ozawa}.

\begin{maintheorem}\label{thm:main}
For every $m\geq2$ and every odd exponent $n\geq557$, the free Burnside group $B(m,n)$ has property $\PHP$. Consequently, its reduced group $C^*$-algebra $\Cr(B(m,n))$, equipped with its canonical trace, is completely selfless, and thus selfless.
\end{maintheorem}

\begin{maincorollary}\label{cor:comparison}
For every $m\geq2$ and every odd exponent $n\geq557$, the algebra $\Cr(B(m,n))$ is simple and uniquely tracial, has stable rank one, and has strict comparison of positive elements by its canonical trace.
\end{maincorollary}

\begin{proof}
Apply \cite[Theorem~3.1]{Robert} to $\Cr(B(m,n))$ with its canonical trace.
\end{proof}

Thus our bound for $C^*$-simplicity and uniqueness of trace matches the best known infinitude threshold. The theorem also recovers the nonamenability of these groups as a consequence of selflessness, as a nontrivial amenable group cannot be $C^*$-simple. 

Beyond free groups, Amrutam, Gao, Kunnawalkam Elayavalli, and Patchell \cite{AGKEP} established selflessness for finitely generated acylindrically hyperbolic groups with trivial finite radical and rapid decay. Ozawa \cite{Ozawa} subsequently removed the rapid decay assumption and obtained complete selflessness in this setting. Vigdorovich obtained examples among higher-rank lattices \cite{VigdorovichLattices}, and subsequently proved selflessness for nontrivial linear groups with trivial amenable radical \cite{VigdorovichLinear}. Bruce and Geffen \cite{BruceGeffen} recently proved that Thompson's group $F$ is $C^*$-selfless, providing an example of this property among groups containing no copy of $\mathbb F_2$. For infinite-dimensional simple unital monotracial $C^*$-algebras, Robert \cite[Theorem~1.1]{RobertComparison} proved that strict comparison by the unique trace implies selflessness; a further characterization appears in \cite{KlissePagliero}.

To our knowledge, the free Burnside groups in Theorem \ref{thm:main} are the first torsion groups whose reduced group $C^*$-algebras are selfless. The exponent law excludes \emph{mixed-identity-freeness} (MIF): the nontrivial word $z^n\in G*\langle z\rangle$, with $\langle z\rangle$ infinite cyclic, vanishes under every substitution $z\mapsto g\in G$. In particular, for finite $m$, these groups are not selfless in the group-theoretic sense of \cite[Definition~3.1]{AGKEP}, which requires a quantitative strengthening of MIF. Ozawa's tensor product theorem already gives non-MIF $C^*$-selfless groups, such as $\mathbb F_2\times\mathbb F_2$ \cite{Ozawa}; the additional feature here is the ordinary group law $x^n=1$.

We also highlight a distinction at the level of boundary actions. Ozawa \cite[Theorem~11]{Ozawa} obtains complete selflessness for infinite countable groups admitting a topologically free minimal extremely proximal action on a compact space. Such an action would force a nonabelian free subgroup by the classical ping-pong argument \cite{Glasner}, and hence cannot exist for
our torsion groups. Their Furstenberg boundary actions are nevertheless topologically free by $C^*$-simplicity \cite{KalantarKennedy}; these actions are strongly proximal, but not extremely proximal.

\subsection*{Proof strategy}

We verify Ozawa's property $\PHP$ by constructing, for each finite test set $F\subseteq G$ and every $N\geq1$, elements $y_1,\ldots,y_N$ and associated subsets satisfying separation and overlap conditions.
Our overlap bound depends only on $F$, independently of $N$, and is therefore stronger than the criterion requires.

To ensure separation, we cover $F^{-1}F\setminus\{1\}$ by finitely many finite cyclic subgroups. A set meeting each right coset of each of these subgroups at most once has pairwise disjoint translates by $F$. We regard each pair $(x,y_i x)$ as a connection and either select one endpoint or discard it. The selected endpoints define the sets in Ozawa's criterion, so we seek a selection satisfying all these coset restrictions.

For a fixed subgroup, a chain alternates between moving within a right coset and traversing a connection in either direction. We require \emph{relative power avoidance} to ensure that chains obeying prescribed cancellation and repetition restrictions cannot revisit a coset. Counting these admissible chains inside a finite family of cosets
restricts the available connections: too many extensions would force a chain visiting more cosets than the family contains. This chain-counting approach builds on \cite{Nonamenability}.

These restrictions yield a decomposition of each coset graph into a forest and a matching. Orienting the forests away from roots gives at most one selected endpoint per coset. We discard the matchings and then further connections to obtain a single selection satisfying all the coset restrictions. The overlap at a point is bounded by the number of connections selecting it plus the number discarded there. The first number is at most one, and the second is bounded in terms of $F$, giving the required estimate independently of $N$.

It remains to construct the same witnesses for all the chosen cyclic subgroups. We use ART's canonical forms to prove a nontriviality lemma for suitable products involving one word and a prescribed finite coefficient set, by exhibiting nonempty canonical representatives. An algebraic reduction then produces the common witnesses as conjugates of one word. The reduction and the canonical-form construction are carried out in Appendix~\ref{sec:canonical}.

\subsection*{Acknowledgments}

I am grateful to Jamie Bell for introducing me to this problem and for his advice throughout the preparation of the paper. I would also like to thank Katrin Tent for her many helpful suggestions and, together with Sam Shepherd, for fruitful discussions on free Burnside groups. I also thank Shirly Geffen for her useful comments.

The author is funded by the Deutsche Forschungsgemeinschaft (DFG, German Research Foundation) under Germany's Excellence Strategy EXC 2044/2 –390685587, Mathematics Münster: Dynamics–Geometry–Structure.

\subsection*{AI-Declaration} 
The mathematical ideas and the development of the proofs in this paper are the author's own, building on previous work in \cite{Nonamenability}. ChatGPT (Astra 6.0) was used to assist with literature searches, proofreading the manuscript, and refining the Appendix. Its contributions primarily concerned the clarity of exposition and the improvement of certain technical bounds. 

\section{Preliminaries and the combinatorial criterion}\label{sec:criterion}

\subsection{Canonical representatives in free Burnside groups}
\label{subsec:canonical-section}

For finite rank $m\geq2$ and odd exponent $n\geq557$, write $\pi:F_m\to B(m,n)$ for the quotient map.
The canonical section $\can:B(m,n)\to F_m$ supplied by \cite[Section~8]{ART} assigns a distinguished reduced word to each group element and satisfies
\[
 \pi\circ\can=\operatorname{id}_{B(m,n)},\qquad \can(1)=1.
\]
For a word $U\in F_m$, we also write $\can(U)$ for $\can(\pi(U))$. Thus
\[
 U=V\text{ in }B(m,n)
 \quad\Longleftrightarrow\quad
 \can(U)=\can(V)\text{ as reduced words}.
\]
In particular, a nonempty canonical representative certifies nontriviality.

ART constructs these representatives through successive quotients of $F_m$, imposing relations $a^n=1$ in increasing construction ranks. Each quotient has its own canonical representatives, and for every fixed word $U$ the resulting sequence eventually stabilizes at $\can(U)$.

For our purposes, the essential features are the estimates controlling periodic subwords under canonical multiplication and the criterion ensuring that a canonical word remains unchanged at later ranks. In Appendix~\ref{sec:canonical}, we use these properties to show that the products used in the power-avoidance argument have nonempty canonical representatives, and hence are nontrivial. The required properties are recalled in Subsection~\ref{subsec:canonical-framework}.

\subsection{Selflessness and strict comparison}

For a discrete group $G$, let $\lambda: G \to \mathcal{U}(\ell^2(G))$ be its left regular representation. The \emph{reduced group $C^*$-algebra} of $G$, denoted by $\Cr(G)$, is the norm closure of the linear span of $\lambda(G)$. Recall that a \emph{state} on $\Cr(G)$ is a positive linear functional $\rho$ such that $\rho(\lambda_{e_G})=1$. A \emph{tracial state} is a state $\tau$ that also satisfies $\tau(ab)=\tau(ba)$ for all $a,b \in \Cr(G)$. The \emph{canonical trace} $\tau_G$ of $G$ is the tracial state on $\Cr(G)$ satisfying
\[
 \tau_G(\lambda_g)=0 \text{ for all } g \in G \setminus \{e_G\}.
\]
A group $G$ is said to be \emph{$C^*$-simple} if $\Cr(G)$ contains no nontrivial closed two-sided ideals. 

Let $A\ne\mathbb C$ be a unital $C^*$-algebra with a faithful tracial state $\tau$, and write
\[
 (B,\sigma)=(A,\tau)*(A,\tau)
\]
for the reduced free product. The pair $(A,\tau)$ is \emph{selfless} if the first-factor embedding
\[
 \iota:(A,\tau)\hookrightarrow(B,\sigma)
\]
is existential \cite[Definition~2.1]{Robert}. Explicitly, this means that, for some ultrafilter $\mathcal U$ on a set $I$, there is a unital trace-preserving $*$-embedding
\[
 j:(B,\sigma)\hookrightarrow
       (A^{\mathcal U},\tau^{\mathcal U})
 \qquad\text{with}\qquad
 j\circ\iota=\Delta_A,
\]
where $\Delta_A$ denotes the diagonal embedding into the norm ultrapower. The index set $I$ need not be countable.

Following \cite[Section~6]{Ozawa}, the pair $(A,\tau)$ is \emph{completely selfless} if $\iota$ is completely existential:
the embedding $j$ can be chosen so that, for every $C^*$-algebra $E$, the map induced by $j$ on elementary tensors extends to a $*$-homomorphism
\[
 B\otimes E\longrightarrow(A\otimes E)^{\mathcal U}.
\]
Here $\otimes$ denotes the minimal tensor product. We say that a group $G$ is \emph{$C^*$-selfless} (respectively, \emph{completely $C^*$-selfless}) if $(\Cr(G),\tau_G)$ is selfless (respectively, completely selfless).

Recall that a unital $C^*$-algebra $A$ has \emph{stable rank one} if its invertible elements are norm-dense in $A$.

We also recall the comparison property appearing in Corollary~\ref{cor:comparison}. Let $A$ be a unital $C^*$-algebra with faithful trace $\tau$, and let $\mathcal K$ denote the compact operators on $\ell^2(\mathbb N)$. For $a,b \in (A \otimes \mathcal K)_+$, we say that $a$ is \emph{Cuntz subequivalent} to $b$, written $a \precsim b$ if there are $v_k \in  A\otimes\mathcal K$ such that $v_k^*bv_k\to a$ in norm. For $c\in(A\otimes\mathcal K)_+$, define the \emph{rank function} \[
 d_\tau(c)=\lim_{k\to\infty}
   (\tau\otimes\operatorname{Tr})(c^{1/k})
 \in[0,\infty],
\]
where $\operatorname{Tr}$ is the usual unnormalized trace. When $A=\mathbb C$, this is the usual rank of a positive compact operator, possibly infinite. We say that $A$ has \emph{strict comparison of positive elements by $\tau$} if, for all $a,b \in (A \otimes \mathcal K)_+$, we have
\[
 d_\tau(a)<d_\tau(b)\quad\Longrightarrow\quad a\precsim b.
\]
The usual definition of strict comparison for simple unital $C^*$-algebras uses all normalized $2$-quasitraces; see \cite[Definition~3.3 and Paragraph~3.4]{APTV}.
For the reduced group $C^*$-algebras in our main theorem, this agrees with the formulation above: selflessness implies simplicity, and the canonical trace is the unique normalized $2$-quasitrace, by \cite[Theorem~3.1]{Robert}.

\subsection{Ozawa's combinatorial property}

\begin{definition}[Property $\PHP$]\label{def:php}
A group $G$ has property $\PHP$ if, for every finite $F\subseteq G$ and $\varepsilon>0$, there are, for every sufficiently large integer $N$, elements $y_1,\ldots,y_N\in G$ and sets $C_i\subseteq D_i\subseteq G$ such that the indexed family
\begin{equation}\label{eq:php-disjoint}
 \{fC_i:f\in F,\ 1\leq i\leq N\}
 \ \cup\
 \{f y_i^{-1}(G\setminus D_i):f\in F,\ 1\leq i\leq N\}
\end{equation}
is pairwise disjoint and
\begin{equation}\label{eq:php-condition}
 \sup_{x\in G}
 \bigl|\{i:x\in D_i\}\cup
        \{i:x\in y_i^{-1}(G\setminus C_i)\}\bigr|
 \leq\varepsilon\sqrt N.
\end{equation}
\end{definition}

The word ``indexed'' means that distinct choices of translate, label $i$, or one of the two displayed families must give disjoint sets. This is the condition in \cite[Section~8]{Ozawa}. The overlap condition says that, at each $x\in G$, all but at
most $\varepsilon\sqrt N$ labels satisfy both
\[
 x\notin D_i \qquad\text{and}\qquad y_i x\in C_i.
\]
Thus, most witnesses move the given point from outside $D_i$ into $C_i$, while \eqref{eq:php-disjoint} separates the prescribed translates.

\begin{theorem}[Ozawa, {\cite[Theorem~14]{Ozawa}}]\label{thm:ozawa}
A group with property $\PHP$ is completely $C^*$-selfless.
\end{theorem}

We will prove a stronger form of the overlap estimate: for each fixed finite set $F$, the following pointwise sum is bounded by a constant depending only on $F$, independently of the number $N$ of witnesses.

\begin{theorem}[Uniform witnesses]\label{thm:uniform}
Let $G=B(m,n)$ satisfy the hypotheses of Theorem~\ref{thm:main}. For every finite $F\subseteq G$ and every integer $N\geq1$, there are $y_1,\ldots,y_N\in G$ and
$C_i\subseteq D_i\subseteq G$ satisfying \eqref{eq:php-disjoint} and
\begin{equation}\label{eq:uniform}
 m(x):=\sum_{i=1}^N\ind_{D_i}(x)
       +\sum_{i=1}^N\ind_{y_i^{-1}(G\setminus C_i)}(x)
 \leq K_F
\end{equation}
for all $x\in G$, where
\[
 K_F=\max\{2,|F|(|F|-1)\}.
\]
\end{theorem}

\begin{proof}[Deduction of Theorem~\ref{thm:main}]
The sum in \eqref{eq:uniform} bounds the cardinality of the union of index sets in \eqref{eq:php-condition}. For fixed $F$ and $\varepsilon>0$, every $N\geq(K_F/\varepsilon)^2$ therefore satisfies Definition~\ref{def:php}. Apply Theorem~\ref{thm:ozawa}.
\end{proof}

\section{Coset connections and encoded chains}
\label{sec:connections}

\subsection{Endpoint choices and cost}\label{subsec:endpoints}

We interpret the sets $C_i\subseteq D_i$ in Ozawa's criterion through choices in a graph with vertex set $G$. This will turn the disjointness condition into a separation
requirement for selected endpoints, and the sum in \eqref{eq:uniform} into a count of contributions at each vertex.

Fix elements $y_1,\ldots,y_N\in G$ with all $2N$ signed elements distinct. For each label $i~\in~\{1,\ldots,N\}$ and $u\in G$, join $u$ to $v=y_i u$ by a \emph{connection}, recorded as the ordered pair $(u,v)$. The group elements $u$ and $v$ are its \emph{endpoints}, and the index $i$ is its \emph{label}, recording which multiplier $y_i$ is used.

For each connection, we may select either endpoint or discard the connection. A connection with a selected endpoint is called \emph{retained} and is oriented towards that endpoint, independently of the order used to record the pair. It contributes one unit of \emph{cost} at its selected endpoint; a discarded connection contributes one unit at each endpoint.

For each label $i$, the sets $C_i\subseteq D_i$ encode these three choices through the position of the second endpoint $v=y_i u$, according to the following table.
\vspace{4mm}
\begin{center}
\begin{tabular}{lll}
\hline
Location of $v=y_i u$ & Decision & Cost\\
\hline
$v\in C_i$ & Select $v$ & One at $v$\\
$v\in D_i\setminus C_i$ & Discard & One at each endpoint\\
$v\notin D_i$ & Select $u$ & One at $u$\\
\hline
\end{tabular}
\end{center}
\vspace{4mm}
Every $v$ is the second endpoint of exactly one connection of label $i$, namely $(y_i^{-1}v,v)$, so the decisions define the sets. The selected second endpoints form $C_i$, and the selected first endpoints form $y_i^{-1}(G\setminus D_i)$: these are exactly the two families appearing in \eqref{eq:php-disjoint}.

At a point $x$, the two connections of label $i$ are $(y_i^{-1}x,x)$ and $(x,y_i x)$. The first contributes at $x$ precisely when $x\in D_i$; the second contributes there precisely when $y_i x\notin C_i$. Hence the sum in \eqref{eq:uniform} is exactly the total cost
at $x$:
\begin{equation}\label{eq:overlap-count}
\begin{aligned}
 m(x)
 &= \#\{\text{retained connections selecting }x\}\\
 &\quad+\#\{\text{discarded connections incident to }x\}.
\end{aligned}
\end{equation}
All $2N$ connections incident to $x$ are distinct by the assumption on the signed elements. A single label may contribute twice to this cost, whereas the union in
\eqref{eq:php-condition} counts that label only once. This explains why the pointwise sum estimate is stronger than Ozawa's overlap condition.

Our task is therefore to make endpoint choices so that the selected second endpoints $C_i$ and the selected first endpoints $y_i^{-1}(G\setminus D_i)$ have pairwise disjoint $F$-translates, as required by \eqref{eq:php-disjoint}, while keeping the total cost at every vertex bounded independently of $N$.

Both requirements involve all labels simultaneously. A bound on the cost contributed by each label separately could still give a total bound growing with $N$.
The arguments below will therefore control the graph containing all connections at once, with a bound depending only on $F$.

\subsection{Connections between right cosets}
\label{subsec:coset-connections}

We now pass from group elements to right cosets. Cosets will help separate the selected endpoints and provide the vertices of the graphs used in the counting argument.

First, let $F\subseteq G$ be finite, and suppose that subgroups $H_1,\ldots,H_r\leq G$ satisfy
\begin{equation}\label{eq:difference-cover}
 F^{-1}F\setminus\{1\}\subseteq\bigcup_{j=1}^r H_j.
\end{equation}
If a set $R\subseteq G$ meets each right $H_j$-coset at most once, then the translates $fR$, $f\in F$, are pairwise disjoint. Indeed, if $fu=f'v$ with $f\ne f'$ and $u,v\in R$, then $u=(f^{-1}f')v$ places $u,v$ in the same right $H_j$-coset for some $j$, which would then give $f=f'$, a contradiction.

To study the choices available in each coset, fix a subgroup $H\leq G$ and elements $Y=\{y_1,\ldots,y_N\}\subseteq G$, where $N\geq1$ and all $2N$ signed elements are distinct. Use the connections $(u,y_i u)$ of Subsection~\ref{subsec:endpoints}, retaining their labels and actual endpoints, and allow traversal in either direction.
Projecting the endpoints to right $H$-cosets gives a graph whose vertices are the cosets and whose edges are the connections. At this stage, loops and multiple connections are allowed. The construction uses individual connections, so no normality assumption on $H$ is required.

A connection $(u,y_i u)$ is \emph{incident to} a right $H$-coset $C$ if at least one of its endpoints belongs to $C$, that is, if $u\in C$ or $y_i u\in C$.
From a current point $x$, any connection incident to its coset $Hx$ can be reached by first moving within that coset to an endpoint $u=hx$, then traversing the connection. Such a step has the form
\[
 x\longmapsto zhx,\qquad z\in Y^{\pm1},\quad h\in H.
\]
The element $h$ is uniquely determined by the chosen endpoint: $h=ux^{-1}$.

A \emph{chain} of length $t\geq0$ consists of a starting point $x_0\in G$ and a symbol list
$\sigma_j=(z_j,h_j)$, $1\leq j\leq t$, determining
\[
 x_j=z_jh_jx_{j-1}.
\]
Its endpoint is
\[
 x_t=z_th_t\cdots z_1h_1x_0.
\]
Thus the list follows the order of successive steps, while the factors appear in reverse order in the resulting product. For $t\geq1$, returning to the initial coset means precisely that
\[
 Hx_t=Hx_0
 \quad\Longleftrightarrow\quad
 z_th_t\cdots z_1h_1\in H.
\]
This is the membership condition that relative power avoidance will exclude for suitable chains.

\subsection{Relative power avoidance}
\label{subsec:avoidance-reduction}
\label{subsec:relative-avoidance}

For an integer $q\geq2$, we say that a word over any alphabet is \emph{$q$-power-free} if it contains no consecutive subword $Q^q$, with $Q$ nonempty. The restriction on repetitions is essential in a group of odd exponent $n$: for $y\ne1$, the word $y^n$ has no adjacent inverse letters but nevertheless represents the identity. Thus the absence of immediate cancellation alone cannot guarantee nontriviality.

\begin{definition}\label{def:relative-avoidance}
Let $H\leq G$, and let $Y=\{y_1,\ldots,y_N\}\subseteq G$ have all $2N$ signed elements distinct. We say that $(H,Y)$ has \emph{relative $q$-power avoidance} if every product
\[
 W=z_th_t\cdots z_1h_1,\qquad
 z_j\in Y^{\pm1},\quad h_j\in H,\quad t\geq1,
\]
lies outside $H$ whenever the following two conditions hold:
\begin{enumerate}
\item Adjacent factors satisfy
\begin{equation}\label{eq:reduced}
 h_{j+1}=1\quad\Longrightarrow\quad z_{j+1}\ne z_j^{-1},
 \qquad1\leq j<t.
\end{equation}
\item The sequence $(z_1,h_1)\cdots(z_t,h_t)$ is $q$-power-free as a word over the alphabet $Y^{\pm1}\times H$.
\end{enumerate}
\end{definition}

In the free product $H*F(Y)$, where $F(Y)$ is free on symbols indexed by $Y$, condition \eqref{eq:reduced} alone would imply $W\notin H$. Relative power avoidance requires this conclusion to persist in $G$ for $q$-power-free symbol lists.

\subsection{Admissible chains}

For the rest of this section, let $G$ be an arbitrary group, $H\leq G$, and $Y=\{y_1,\ldots,y_N\}\subseteq G$, where $N\geq1$, the $2N$ signed elements are distinct, and $(H,Y)$ has relative $15$-power avoidance.
 
\begin{definition} \label{def:admissible-chains}
    A chain with steps
    \[ x_k=z_kh_kx_{k-1},\qquad
 z_k\in Y^{\pm1},\quad h_k\in H,\quad 1\leq k\leq t,\]
 is \emph{admissible} if:
 \begin{enumerate}
\item it has no immediate reversal, meaning that
\[
 h_{k+1}=1\quad\Longrightarrow\quad z_{k+1}\ne z_k^{-1},
 \qquad 1\leq k<t;
\]
\item its symbol list $(z_1,h_1)\cdots(z_t,h_t)$ is $15$-power-free.
\end{enumerate}
\end{definition}

The first condition is precisely \eqref{eq:reduced}. Indeed, reversing the connection just traversed means starting at the point just reached, without an intervening $H$-move, and using the inverse letter.

Every contiguous subchain of an admissible chain is admissible: the reduction condition is inherited by adjacent steps, and a contiguous subword of a $15$-power-free list is again $15$-power-free.

\begin{lemma}\label{lem:no-return}
An admissible chain of length $t$ visits $t+1$ distinct right $H$-cosets.
\end{lemma}

\begin{proof}
Suppose that $Hx_i=Hx_j$ for some $0\leq i<j\leq t$. The steps between $x_i$ and $x_j$ give
\[
 x_j=Wx_i,\qquad
 W=z_jh_j\cdots z_{i+1}h_{i+1}.
\]
Since $Hx_i=Hx_j$, there is $h\in H$ such that $x_j=hx_i$. Cancelling $x_i$ on the right gives $W=h\in H$.

On the other hand, the subchain from $x_i$ to $x_j$ is nonempty and admissible. Its symbol list therefore satisfies both conditions in the definition of relative $15$-power avoidance, which gives $W\notin H$ and thus the desired contradiction.
\end{proof}

\begin{lemma}\label{lem:simple}
No connection has both endpoints in one right $H$-coset, and at most one connection joins any two distinct right $H$-cosets. Consequently, for any family of $k$ right $H$-cosets, at most $\binom{k}{2}$ connections have both endpoints in their union.
\end{lemma}

\begin{proof}
A connection inside one coset would give an admissible one-step return, contradicting Lemma~\ref{lem:no-return}.

Suppose that two distinct connections join cosets $C$ and $D$. Traverse the first from $C$ to $D$. From the endpoint reached, move within $D$ by an element of $H$ to the endpoint of the second connection, then traverse that connection back to $C$.
This gives a two-step chain without immediate reversal, since the two connections are distinct. Its symbol list is automatically $15$-power-free, so the chain is admissible.
The return to $C$ contradicts Lemma~\ref{lem:no-return}.
\end{proof}

We have excluded returns for admissible chains. The counting argument below uses this restriction to control the number and arrangement of connections inside finite families of cosets.

\subsection{The counting argument}

Let $\mathcal C$ be a nonempty finite family of right $H$-cosets and put $\Omega_{\mathcal C}=\bigcup_{C\in\mathcal C}C$. A connection is \emph{internal} if both endpoints belong to $\Omega_{\mathcal C}$.

Consider the graph whose vertices are the cosets in $\mathcal C$ and whose edges are the internal connections: a connection $(u,y_i u)$ joins the vertices $Hu$ and $Hy_i u$.
By Lemma~\ref{lem:simple}, this graph is finite and simple. For $C\in\mathcal C$, write
\begin{equation}\label{eq:degree}
 d_{\mathcal C}(C)
 =\#\{\text{internal connections incident to }C\}
\end{equation}
for its degree.

For $x\in C$, an incident internal connection with endpoint $u\in C$ determines a unique step $x\mapsto zhx$: take $h=ux^{-1}\in H$, with $z$ determined by the label and traversal direction.
Before imposing admissibility, this correspondence gives
\begin{equation}\label{eq:available}
 \#\{(z,h)\in Y^{\pm1}\times H:
       zhx\in\Omega_{\mathcal C}\}
 =d_{\mathcal C}(C).
\end{equation}
After arrival along an internal connection, excluding its immediate reversal leaves exactly $d_{\mathcal C}(C)-1$ choices.

By Lemma~\ref{lem:no-return}, an admissible chain inside $\Omega_{\mathcal C}$ has length at most $|\mathcal C|-1$. We use this bound to obtain the following alternative.

\begin{proposition}\label{prop:small}
Every nonempty finite family $\mathcal C$ of right $H$-cosets contains either a coset with at most one incident internal connection, or two connected cosets each having exactly two  incident internal connections.
\end{proposition}

\begin{proof}
Suppose otherwise.  We will show that enough two-step extensions survive the exclusion of forbidden repetitions to produce arbitrarily long admissible chains inside $\Omega_{\mathcal C}$, contradicting
Lemma~\ref{lem:no-return}. All coset degrees are at least two, and degree-two cosets are nonadjacent.
By \eqref{eq:available}, every chain without immediate reversal, including the empty chain, has at least two two-step extensions without immediate reversal inside $\Omega_{\mathcal C}$: each first step has a continuation; if there is only one available first step, the current coset has degree two, so the next coset has degree at least three and provides at least two continuations.

We therefore group successive symbols into pairs $\eta_j=(\sigma_{2j-1},\sigma_{2j})$.
Fix $x_0\in\Omega_{\mathcal C}$, and let $\alpha_t$ count chains of length $2t$ from $x_0$, entirely inside $\Omega_{\mathcal C}$, with no immediate reversal and a $5$-power-free list of pairs.
These counts are finite by \eqref{eq:available}, and $\alpha_0=1$.

Every counted chain is admissible. Indeed, a factor $P^{15}$, with $p=|P|$, contains five
equal pieces of length $2p$ aligned with the pairing: shift its starting position by at most one symbol and use $15p-1\geq10p$. These pieces would give a fifth power in the paired list.

There are at least $2\alpha_t$ candidate extensions. Since each appends one pair, any newly created fifth power must be a suffix, giving a list $UQ^5$. For fixed $q=|Q|$, measured in pairs, delete the final four copies of $Q$. The remaining prefix $UQ$ is a prefix of the original counted chain, so it is counted by $\alpha_{t+1-4q}$. Its last $q$ pairs recover $Q$ and hence the entire failed extension $UQ^5$. Thus, the deletion is injective for each fixed $q$.
A union bound therefore gives
\begin{equation}\label{eq:count}
 \alpha_{t+1}\geq2\alpha_t-
 \sum_{q=1}^{\lfloor(t+1)/5\rfloor}\alpha_{t+1-4q}.
\end{equation}

Starting from $\alpha_1\geq2$, induction yields $\alpha_{t+1}\geq\tfrac32\alpha_t$.
Indeed, the previous growth estimates imply
\[
 \sum_{q=1}^{\lfloor(t+1)/5\rfloor}\alpha_{t+1-4q}
 \leq
 \alpha_t\sum_{q\geq1}(2/3)^{4q-1}
 =\frac{24}{65}\alpha_t
 <\frac12\alpha_t,
\]
which proves the induction step through \eqref{eq:count}.

Thus $\alpha_t>0$ for every $t$, giving arbitrarily long admissible chains inside $\Omega_{\mathcal C}$.
This contradicts Lemma~\ref{lem:no-return}, which bounds their length by $|\mathcal C|-1$.
\end{proof}

\begin{remark}[Nonamenability]\label{rem:nonamenability}

For $G=B(m,n)$ as in Theorem~\ref{thm:main}, choose
$Y=\{y_1,y_2\}$ using Theorem~\ref{thm:common-avoidance}
with $H_1=\{1\}$ and $N=2$. For finite $A\subseteq G$,
let $e(A)$ count the connections with both endpoints in $A$.
Deleting the configurations in Proposition~\ref{prop:small}
gives $e(A)\leq\tfrac32|A|$ by induction on $|A|$. Hence
\[
 |y_1A\setminus A|+|y_2A\setminus A|
 =2|A|-e(A)\geq\tfrac12|A|.
\]
By F\o lner's criterion, $G$ is nonamenable, recovering
the result of \cite{Nonamenability}.

\end{remark}

\section{Construction of Ozawa's sets}\label{sec:forests}

We now use Proposition~\ref{prop:small} to construct the endpoint choices required by Ozawa's criterion. For a single subgroup, we decompose the coset graph into a forest and a matching. Discarding the matching and orienting the forest gives at most one selected endpoint in each coset.
We then construct endpoint choices such that, for every prescribed subgroup $H_j$, each right $H_j$-coset contains the selected endpoint of at most one retained connection.
These choices determine the sets $C_i\subseteq D_i$ through the correspondence described in Subsection~\ref{subsec:endpoints}.

\subsection{Decomposition of the coset graph}
\label{subsec:forest-decomposition}

A \emph{forest} is a graph with no cycles, and a \emph{matching} is a set of edges with pairwise disjoint endpoints. Orienting each component of a forest away from a chosen root gives at most one incoming edge at each vertex. Thus, in a coset graph, this orientation selects at most one endpoint in each coset. A discarded matching contributes at most one discarded connection at each coset.

The two configurations found in Proposition~\ref{prop:small} allow us to construct such a decomposition by deleting and then restoring one or two vertices.
\begin{proposition}
\label{prop:forest}
Let $G$ be a group, let $H\leq G$, and let $Y=\{y_1,\ldots,y_N\}\subseteq G$, with $N\geq1$ and all
$2N$ signed elements distinct. Suppose that $(H,Y)$ has relative $15$-power avoidance. Then the right-$H$-coset connection graph is simple, and its edge set $E$ admits a partition
\[
 E=\mathcal T\sqcup\mathcal M,
\]
where $\mathcal T$ is a forest and $\mathcal M$ is a matching.
\end{proposition}

\begin{proof}
The graph is simple by Lemma~\ref{lem:simple}. First consider a finite induced subgraph and argue by
induction on its number of vertices. The empty graph requires no choices.

By Proposition~\ref{prop:small}, every nonempty such graph has either a vertex of degree at most one or two adjacent vertices of degree two. In the first case, remove the vertex and decompose the remaining graph by induction. Restore the removed vertex and place its edge, if any, in the forest. This adds an isolated vertex or a leaf, so the decomposition is preserved.

In the second case, let $C,D$ be adjacent degree-two vertices, with other incident edges $AC$ and $DB$. Here $A,B$ lie outside $\{C,D\}$ and may coincide. Remove $C,D$ and decompose the remaining graph by induction. Restore $AC$ and $DB$ to the forest, making $C$ and $D$ leaves, and place $CD$ in the matching. Since both endpoints of $CD$ were absent from the smaller graph, the matching condition is preserved.

For the full graph, consider the compact space $\{\mathrm{T},\mathrm{M}\}^{E}$ of edge assignments. Require that no cycle be assigned entirely to $\mathrm{T}$ and that no two incident edges both be assigned to $\mathrm{M}$. These are closed conditions, each involving finitely many edges, and any finite collection of conditions involves only finitely many vertices. Decompose the induced subgraph on those vertices by the finite case and extend the assignment arbitrarily to the remaining edges. Thus, the conditions have the finite intersection property, and compactness
gives the required decomposition.

\end{proof}

\subsection{Simultaneous endpoint choices}
\label{subsec:simultaneous-endpoints}

For each $H_j$, a connection $(x,y_i x)$ joins the cosets $H_jx$ and $H_jy_i x$, so all coset graphs have the same edge set $E$. We seek a single endpoint selection such that each right $H_j$-coset contains the selected endpoint of at most one retained connection, while bounding the number of discarded connections incident to each point of $G$.

After discarding the matchings, we combine the forest orientations using one observation: orienting a tree away from a root selects every vertex except that root.
The following lemma makes this precise.
\begin{lemma}
\label{lem:common}
Let $G$ be a group and let $Y=\{y_1,\ldots,y_N\}\subseteq G$, with $N\geq1$ and all $2N$ signed elements distinct.
Let
\[
 E=\{(x,y_i x):x\in G,\ 1\leq i\leq N\}
\]
be the common set of connections.
Suppose that $H_1,\ldots,H_r\leq G$, with $r\geq1$, have the following property: for each $j$, the right-$H_j$-coset connection graph is simple and admits a decomposition
\[
 E=\mathcal T_j\sqcup\mathcal M_j
\]
into a forest and a matching. Then, some connections can be discarded and an endpoint selected for each remaining connection, so that:
\begin{enumerate}
\item for every $j$, each right $H_j$-coset contains the selected endpoint of at most one retained connection;
\item at most $2r-1$ discarded connections are incident to any given point of $G$.
\end{enumerate}
\end{lemma}

\begin{proof}
First discard all connections in
\[
 \mathcal M=\bigcup_{j=1}^r\mathcal M_j,
 \qquad E_0=E\setminus\mathcal M.
\]
At most $r$ discarded connections are incident to any point $x\in G$, since each matching $\mathcal M_j$ has at most one edge incident to the coset $H_jx$.
For every $j$, the remaining coset graph is a forest, because $E_0\subseteq\mathcal T_j$.

For each $j$, root every component of this coset forest and orient its edges away from the root. Orient each actual connection towards its endpoint in the terminal coset, and let $R_j$ be the set of selected endpoints. At most one edge enters each coset, so each right $H_j$-coset contains the selected endpoint of at most one connection. In particular, each point of $G$ has at most one incoming edge.

We claim that, in each connected component $T$ of the graph $(G,E_0)$, this orientation selects every vertex except one. To see this, start at any vertex and follow incoming edges backwards for as long as possible. At each step, the image in the coset forest moves one edge towards its chosen root. The process therefore stops after finitely many steps, at a vertex with no incoming edge.

Moreover, adjacent vertices have the same terminal vertex for this process: if an edge points from $u$ to $v$, the backward path from $v$ first visits $u$. By connectedness, all vertices of $T$ therefore lead to the same vertex, which we denote by $v_j(T)$. This is the unique vertex of $T$ with no incoming edge, and hence
\[
 R_j\cap V(T)=V(T)\setminus\{v_j(T)\}.
\]

Now, keep the orientation corresponding to $j=1$. In each component $T$, discard the unique incoming
edge at each distinct vertex among $v_2(T),\ldots,v_r(T)$ other than $v_1(T)$.
Leave all remaining orientations unchanged. The resulting set $R$ of selected endpoints satisfies
\[
 R\cap V(T)
 =V(T)\setminus\{v_1(T),\ldots,v_r(T)\}.
\]
Thus $R\subseteq R_j$ for every $j$, so $R$ meets each right $H_j$-coset at most once. Since deleting edges from the first orientation preserves the property that each point is selected by at most one connection, this proves the first assertion.

At most $r-1$ further edges were discarded in each original component $T$ of $(G,E_0)$. Together with
the at most $r$ incident connections discarded with the matchings, this gives at most $2r-1$ discarded connections incident to any point of $G$.
\end{proof}

\subsection{Witnesses from relative avoidance}
\label{sec:witnesses}

We now use the simultaneous endpoint choices to construct the sets required by Ozawa's criterion.

\begin{proposition}[Witnesses from relative avoidance]
\label{prop:witnesses}
Let $G$ be a group, let $F\subseteq G$ be finite, and let $N\geq1$. Suppose that subgroups $H_1,\ldots,H_r\leq G$, with $r\geq1$, satisfy
\[
 F^{-1}F\setminus\{1\}\subseteq\bigcup_{j=1}^rH_j,
\]
and that $Y=\{y_1,\ldots,y_N\}\subseteq G$ has all $2N$ signed elements distinct, with every $(H_j,Y)$ satisfying relative $15$-power avoidance. Then there are sets $C_i\subseteq D_i\subseteq G$ satisfying \eqref{eq:php-disjoint} and $m(x)\leq2r$ for all $x \in G$.

\end{proposition}

\begin{proof}
Apply Proposition~\ref{prop:forest} to each $H_j$, using the same connections
\[
 E=\{(x,y_i x):x\in G,\ 1\leq i\leq N\}.
\]
Lemma~\ref{lem:common} gives discarded connections and endpoint choices for the retained connections. Let $R$ be the set of selected endpoints. No point is selected twice, and $R$ meets each right $H_j$-coset at most once. By the assumption on $F^{-1}F \setminus \{1\}$ and the observation in Subsection~\ref{subsec:coset-connections}, the translates $fR$, $f\in F$, are pairwise disjoint.

Record the endpoint choices as in Subsection~\ref{subsec:endpoints}:
\begin{align*}
 C_i&=\{y_i x:(x,y_i x)\text{ is retained and selects }y_i x\},\\
 D_i&=C_i\cup\{y_i x:(x,y_i x)\text{ is discarded}\}.
\end{align*}
The selected first endpoints of label $i$ form $y_i^{-1}(G\setminus D_i)$.
Since no point is selected twice, these sets and the $C_i$ partition $R$ as an indexed family. Their indexed $F$-translates are therefore pairwise disjoint, giving \eqref{eq:php-disjoint}.

Finally, at most one retained connection selects any given point, and at most $2r-1$ discarded connections meet that point. The cost identity \eqref{eq:overlap-count} therefore gives
\[
 m(x)\leq1+(2r-1)=2r
\]
for all $x \in G$.
\end{proof}

\section{Application to free Burnside groups}\label{sec:application}

To apply Proposition~\ref{prop:witnesses}, we need common relative-avoidance witnesses for a finite family of finite subgroups. The following theorem provides them.

\begin{theorem}[Common relative avoidance]
\label{thm:common-avoidance}
Let $G=B(m,n)$, where $m\geq2$ and $n\geq557$ is odd. Let $H_1,\ldots,H_r\leq G$ be finite subgroups, with $r\geq1$.
For every integer $N\geq1$, there are pairwise conjugate elements $y_1,\ldots,y_N\in G$, with all $2N$ signed elements distinct, such that $(H_j,Y)$ has relative $15$-power avoidance for every $j$, where $Y=\{y_1,\ldots,y_N\}$.
\end{theorem}

The theorem is proved in Appendix~\ref{sec:canonical}, using the canonical forms recalled in
Subsection~\ref{subsec:canonical-section}. We now combine it with Proposition~\ref{prop:witnesses} to obtain the uniform overlap bound.

\begin{proof}[Proof of Theorem~\ref{thm:uniform}]
Let $G=B(m,n)$ satisfy the hypotheses of the theorem, and fix a finite set $F=\{f_1,\ldots,f_k\}\subseteq G$.

For $k\geq2$, let $H_1,\ldots,H_r$ be the distinct subgroups among
\[
 \langle f_a^{-1}f_b\rangle,
 \qquad 1\leq a<b\leq k.
\]
These subgroups are finite because every element of $G$ has order dividing $n$. Each contains both
$f_a^{-1}f_b$ and its inverse $f_b^{-1}f_a$, so
\[
 F^{-1}F\setminus\{1\}
 \subseteq\bigcup_{j=1}^r H_j,
 \qquad r\leq\binom{k}{2}.
\]
For $k\leq1$, take $r=1$ and $H_1=\{1\}$; the same covering condition holds because $F^{-1}F\setminus\{1\}$ is empty.

With these subgroups fixed, let $N\geq1$. Theorem~\ref{thm:common-avoidance} supplies
$Y=\{y_1,\ldots,y_N\}$ with all $2N$ signed elements distinct and such that every $(H_j,Y)$ has relative $15$-power avoidance. Proposition~\ref{prop:witnesses} therefore gives sets $C_i\subseteq D_i\subseteq G$ satisfying \eqref{eq:php-disjoint} and, for every $x\in G$,
\[
 m(x)
 \leq 2r
 \leq \max\{2,k(k-1)\}
 =K_F.
\]
This proves \eqref{eq:uniform} for every $N\geq1$.
\end{proof}

\appendix
\raggedbottom

\section{Relative power avoidance via canonical forms} \label{sec:canonical}

Throughout this appendix, let $G=B(m,n)$, where $m\geq2$ and $n\geq557$ is odd, and write $\tau=15$ for the nesting constant in \cite{ART}.

Our goal is to prove Theorem~\ref{thm:common-avoidance}. We first reduce the construction of common witnesses to a nontriviality statement for products involving
one word and coefficients from a prescribed finite set (Lemma~\ref{lem:one-word-avoidance}). The remainder of the appendix proves this lemma using ART's canonical forms. We construct the word so that nonempty periodic subwords survive in the canonical representatives of the required products, thereby certifying their nontriviality.

\subsection{From one word to common relative witnesses}
\label{subsec:one-word-reduction}

We construct the common witnesses as conjugates $y_i=t_iwt_i^{-1}$. The lemma below is used first to choose the conjugators $t_i$, and then to choose $w$ for the finite
coefficient set arising when these conjugates are expanded.

The lemma asks a single element $w$ to behave like a new free letter in products with coefficients from a prescribed finite set $E$. We consider products with no immediate
cancellation and exclude long repetitions in the record of signs and intervening coefficients. Once $E$ is fixed, the same $w$ makes every such product nontrivial, with no prescribed bound on its length.

\begin{lemma}
\label{lem:one-word-avoidance}
Let $E\subseteq G$ be finite and contain $1$. There exists $w\in G$ with the following property. Let $t\geq1$, $e_0,\ldots,e_t\in E$, and $\varepsilon_1,\ldots,\varepsilon_t\in\{\pm1\}$.
Suppose that
\[
 e_j=1\quad\Longrightarrow\quad
 \varepsilon_{j+1}=\varepsilon_j
\]
for $1\leq j<t$, and that the interior record
\[
 (\varepsilon_1,e_1)\cdots
 (\varepsilon_{t-1},e_{t-1})
\]
is $\tau$-power-free. Then
\[
 e_0w^{\varepsilon_1}e_1\cdots
 w^{\varepsilon_t}e_t\ne1.
\]
For $t=1$, the interior record is empty.
\end{lemma}

The repetition restriction is essential: the relation $w^n=1$ is excluded because its interior record $(+1,1)^{n-1}$ contains a $\tau$-th power.
We prove this lemma in Subsection~\ref{subsec:one-word-proof}. The first and last coefficients are left out of the record. They can therefore be changed within $E$ without affecting the power-freeness condition.

For the next lemma, $\langle x\rangle$ is infinite cyclic. We count occurrences of $x^{\pm1}$ in reduced free-product forms by expanding each power $x^k$ into $|k|$ letters.

\begin{lemma}
\label{lem:short-substitutions}
Let $P\subseteq G*\langle x\rangle$ be finite. Suppose that the reduced form of every element of $P$ contains at most two occurrences of $x^{\pm1}$. Then there exists $u\in G$ such that evaluation at $x=u$, fixing $G$, is injective on $P$.
\end{lemma}

\begin{proof}
For distinct $p,q\in P$, reduce $pq^{-1}$ in $G*\langle x\rangle$. A nonidentity constant remains
nonidentity under every evaluation. Otherwise, the reduced expression has the form
\[
 a_0x^{\varepsilon_1}a_1\cdots
 x^{\varepsilon_\ell}a_\ell,
\]
with $1\leq\ell\leq4$. Reducedness gives the reduction condition of Lemma~\ref{lem:one-word-avoidance}.
Collect the coefficients of all these expressions in one finite set $E$, including $1$. Their interior records have length at most three and are therefore $\tau$-power-free. Lemma~\ref{lem:one-word-avoidance} supplies one evaluation under which none of the products $pq^{-1}$ becomes trivial.
\end{proof}

Expanding a product in the conjugates replaces a coefficient $h$ between $y_i^\varepsilon$ and $y_j^\delta$ by $t_i^{-1}ht_j$. We choose the conjugators so that two different triples $(i,h,j)$ cannot give the same nonidentity element. This will allow us to compare repetitions in the expanded product with repetitions in the original symbol list.

\begin{lemma}
\label{lem:separated-conjugators}
Let $D\subseteq G$ be finite and contain $1$, and let $N\geq1$. There exist $t_1=1,t_2,\ldots,t_N\in G$ such that, for $d\in D$ and $1\leq i,j\leq N$,
\[
 t_i^{-1}dt_j=1
 \quad\Longleftrightarrow\quad
 i=j\text{ and }d=1.
\]
Moreover, every nonidentity value $t_i^{-1}dt_j$ uniquely determines $i$, $d$, and $j$.
\end{lemma}

\begin{proof}
We reason by induction on the number of conjugators, starting with $t_1=1$.
Suppose that $t_1,\ldots,t_k$ have been chosen, and temporarily put $t_{k+1}=x$ in $G*\langle x\rangle$. Apart from the identity expression $x^{-1}1x=1$, the new mixed coefficients have the forms
\[
 (t_i^{-1}d)x,\qquad
 x^{-1}(dt_j),\qquad
 x^{-1}dx.
\]
Here $i,j\leq k$ and $d\in D$, with $d\ne1$ in the third family. Free-product normal forms distinguish these families from one another and from the old mixed coefficients, which lie in $G$.
Within the first family, equality reduces to $t_i^{-1}d=t_{i'}^{-1}e$; the inductive hypothesis, using $t_1=1$ as the right index, gives $i=i'$ and $d=e$.
The second family is treated with $t_1=1$ as the left index. In the third family, cancellation gives $d=e$. Thus the required separation holds with the formal conjugator $x$.

Let $P$ be the finite set of all mixed-coefficient values for these $k+1$ conjugators in the free product. Every element of $P$ contains at most two occurrences of $x^{\pm1}$ in reduced form. Lemma~\ref{lem:short-substitutions} provides an evaluation $x\mapsto t_{k+1}\in G$ preserving all distinct values in $P$. This completes the induction.
\end{proof}

\begin{proof}[Proof of Theorem~\ref{thm:common-avoidance}]
Put $D=\bigcup_{j=1}^rH_j$ and choose $t_1=1,t_2,\ldots,t_N$ by Lemma~\ref{lem:separated-conjugators}. The set
\[
 E=\bigcup_{i,k=1}^N t_i^{-1}Dt_k
\]
is finite and contains $1$. Apply Lemma~\ref{lem:one-word-avoidance} to obtain $w$, and
set $y_i=t_iwt_i^{-1}$ for $1\leq i\leq N$. Since $t_1=1$, the set $E$ also contains $Dt_i$ and $t_i^{-1}D$ for every $i$.

The $2N$ signed witnesses are distinct. Indeed, take signs $\varepsilon,\delta\in\{\pm1\}$
and indices $i,k$ such that $i\ne k$ or $\varepsilon\ne\delta$. Put $c=t_i^{-1}t_k$.
Then $c,c^{-1}\in E$, and
\[
 t_i^{-1}y_i^\varepsilon(y_k^\delta)^{-1}t_i
 =w^\varepsilon c\,w^{-\delta}c^{-1}.
\]
If $c=1$, Lemma~\ref{lem:separated-conjugators} gives $i=k$, and hence $\varepsilon=-\delta$. Thus the displayed product satisfies the reduction condition of Lemma~\ref{lem:one-word-avoidance}. Its interior record has length one, so that lemma shows that the product is nonidentity.

Fix $H=H_j$, and let
\[
 (y_{i_1}^{\varepsilon_1},h_1)\cdots
 (y_{i_t}^{\varepsilon_t},h_t)
\]
be a reduced, $\tau$-power-free relative symbol list, where $t\geq1$ and $h_1,\ldots,h_t\in H$.
Put
\[
 W=y_{i_t}^{\varepsilon_t}h_t\cdots
   y_{i_1}^{\varepsilon_1}h_1.
\]
For $2\leq k\leq t$, define
\[
 c_k=t_{i_k}^{-1}h_kt_{i_{k-1}}.
\]
We claim that the record
\[
 (\varepsilon_2,c_2)\cdots(\varepsilon_t,c_t)
\]
is also $\tau$-power-free, with the record understood to be empty when $t=1$.

Suppose that a nonempty block repeats $\tau$ times in this record. Each nonidentity coefficient $c_k$ determines $i_k,h_k,i_{k-1}$ by Lemma~\ref{lem:separated-conjugators}, whereas $c_k=1$ forces $h_k=1$ and $i_k=i_{k-1}$. If the repeated block contains a nonidentity coefficient, these data determine all indices and $H$-coefficients in each copy: the identity coefficients propagate the determined indices through the intervening strings, including those at either end of the block. The copies therefore give identical blocks of original symbols. If every coefficient in the block is identity, all indices throughout the repeated segment coincide and all its $H$-coefficients are identity.
In either case, the signs already repeat, so the original symbol list contains a $\tau$-th power. This contradiction proves the claim.

For any $h_0\in H$, expanding the conjugates gives
\[
 h_0W=(h_0t_{i_t})w^{\varepsilon_t}c_t
 \cdots w^{\varepsilon_2}c_2
 w^{\varepsilon_1}(t_{i_1}^{-1}h_1),
\]
with no intermediate coefficients when $t=1$. All coefficients belong to $E$. If $c_k=1$, then $h_k=1$ and $i_k=i_{k-1}$, so reducedness of the original list gives $\varepsilon_k=\varepsilon_{k-1}$. The interior record of this expansion is the reversal of the record just considered, and is therefore $\tau$-power-free.
Lemma~\ref{lem:one-word-avoidance} gives $h_0W\ne1$. Since $h_0\in H$ was arbitrary, $W\notin H$. This proves relative $\tau$-power avoidance for every $H_j$, completing the proof.
\end{proof}

\subsection{Canonical forms and periodic estimates}
\label{subsec:canonical-framework}

We recall the canonical-form facts needed to prove Lemma \ref{lem:one-word-avoidance}.
We use ART's construction on a finite alphabet, which suffices here because the lemma involves only a finite coefficient set. Indeed, we can choose finitely many basis generators, say $k\geq2$, whose generated subgroup contains all these coefficients. The induced map $B(k,n)\to G$ is injective, since sending the remaining basis generators to the identity defines a left inverse. Thus, nontriviality proved in this basis subgroup persists in $G$. For the remainder of the appendix, we work with $G=B(m,n)$, where $2\leq m<\infty$.

ART organizes the defining relations into sets $\mathrm{Rel}_r$ and forms the successive quotients
\[
 G_r=F_m/\langle\!\langle
       \mathrm{Rel}_1\cup\cdots\cup\mathrm{Rel}_r
       \rangle\!\rangle,
\]
with $G_0=F_m$. Each relator in $\mathrm{Rel}_r$ has the form $a^n$, where the \emph{period} $a$ is nonempty, cyclically reduced, and primitive, where the latter here means that $a$ is not a proper power. The relator sets are closed under cyclic shifts and inverses. Their construction ranks record the nesting of periodic subwords and are unrelated to the number of generators.

For each word $U$, let $\can_r(U)$ be its rank-$r$ canonical representative, and write $\mathrm{Can}_r$ for the set of these representatives. The sequence $\can_r(U)$ eventually stabilizes and gives the canonical section $\can:G\to F_m$
\cite[Proposition~8.2, Corollary~8.7, and Theorem~8.1]{ART}.
The sets $\mathrm{Can}_r$ are nested by IH~2. Consequently, every final canonical representative belongs to $\mathrm{Can}_r$ for every $r\geq0$.
For $s\geq1$, a nonempty primitive cyclically reduced word $a$ is \emph{cyclically canonical at rank $s$} if $a^n\notin\mathrm{Rel}_1\cup\cdots\cup\mathrm{Rel}_s$ and $a^\tau$ occurs as a subword of some word in $\mathrm{Can}_s$ \cite[Definition~4.2]{ART}.
Factorizations of reduced words are understood without cancellation, and $\red$ denotes free reduction. In $G_s$, coefficients $e\in G$ are represented by $\can(e)$, while constructed words are represented by their displayed expressions.
We retain the same symbols for these lifts and specify the quotient in which equalities hold.

A \emph{fractional power} of $a$ is a subword of a power of $a$. For an occurrence $u$ with this specified period, put
\[
 \Lambda_a(u)=\frac{|u|}{|a|}.
\]
For a reduced word $U$, let $\ell_a(U)$ be the supremum of $|u|/|a|$ over all fractional-power occurrences $u$ in $U$ with period a cyclic shift of $a^{\pm1}$.
Let $\ell_r(U)$ be the supremum of $\ell_a(U)$ over all periods $a$ of rank $r$. We include occurrences shorter than one full period, with their period specified, and set empty suprema equal to zero.
Both functions are subadditive under concatenation without cancellation. Indeed, cutting a periodic occurrence at a junction gives pieces with the same period, each counted by our definition even if shorter than one full period.

An $a$-periodic occurrence is \emph{maximal} if it cannot be extended in its chosen periodic phase. An adjacent letter preventing extension is a \emph{stopping letter}.
For an occurrence of at least one full period, its sign records whether its phase comes from $a$ or $a^{-1}$. These alternatives are distinct because a nonempty cyclically reduced word is not a cyclic shift of its inverse.

We use the following four properties. All IH numbers refer to Section~3 of \cite{ART}.

\begin{enumerate}
\item\label{fact:canonical-equality}
\emph{Equality and survival.}
Canonical forms identify equality in $G_r$ and commute
with inversion, by IH 7--8.
Set
\begin{equation}\label{eq:canonical-threshold}
 T=\frac n2-5\tau-2>12\tau+21.
\end{equation}
By IH~4, a word $U\in\mathrm{Can}_{r-1}$ with $\ell_r(U)\leq T$ already belongs to $\mathrm{Can}_r$.
Thus, if a nonempty word $U\in\mathrm{Can}_s$ satisfies $\ell_r(U)\leq T$ for every $r>s$, it remains unchanged at every later rank and represents a nonidentity element of $G$.

\item\label{fact:canonical-triangle}
\emph{Canonical multiplication.}
For $U,V\in\mathrm{Can}_s$, IH~10 gives literal factorizations
\begin{equation}\label{eq:canonical-triangle}
 U=U'D_1X,\qquad V=X^{-1}D_2V',
 \qquad \can_s(UV)=U'D_3V'.
\end{equation}
Here $XX^{-1}$ is the maximal free cancellation in the product $UV$. The words $D_1,D_2,D_3$ form a \emph{canonical triangle}: they satisfy
$D_1D_2=D_3$ in $G_s$ and
\[
 \ell_i(D_j)<\tau
\]
for every $i>s$ and $j\in\{1,2,3\}$. The analogous prefix and suffix formulas use a canonical triangle at the exposed end, by IH~12.
These estimates control periodic measure: triangle sides can have large ordinary word length.

\item\label{fact:canonical-locality}
\emph{Changing periodic blocks.}
Let $a$ be a period of rank greater than $s$. If
\[
 L_1a^\tau R_1,\ L_2a^\tau R_2\in\mathrm{Can}_s,
\]
then
\[
 L_1a^qR_2\in\mathrm{Can}_s
\]
for every integer $q\geq\tau$, by IH 11.

A surviving block $a^\tau$ also separates the effects of changes in the two exterior contexts. More precisely, in \eqref{eq:canonical-triangle},
suppose that $V'=Ma^\tau R$. If
\[
 V_1=X^{-1}D_2Ma^\tau R_1\in\mathrm{Can}_s,
\]
then
\[
 \can_s(UV_1)=U'D_3Ma^\tau R_1.
\]
The corresponding statement at the left end follows by inversion
\cite[Corollary~3.9]{ART}.

\item\label{fact:period-overlap}
\emph{Small cancellation condition.}
By IH~5, incompatible periodic occurrences of the same rank overlap by less than $\tau+1$ periods, measured with respect to either period. Here incompatible means that their bi-infinite periodic continuations, aligned at the overlap, are different. Occurrences of different ranks overlap by less than two periods of the higher rank.
For a fixed primitive period $a$, two distinct maximal $a$-periodic occurrences overlap by less than $|a|$: otherwise a full common period would identify their phases, and their union would be one periodic occurrence.
\end{enumerate}

\subsection{Finite coefficients and controlled multiplication}

We choose a stage where all the coefficient words already represent elements of order dividing $n$. These words have fixed lengths, so we can then move further along the construction until every later period is long compared with all of them.

\begin{lemma}
\label{lem:finite-stage-data}
Let $E\subseteq G$ be finite and let $\eta>0$. There exists an integer $R\geq1$ with the following
properties:
\begin{enumerate}
\item For every $e\in E$, the word $\can(e)$ represents an element of order dividing $n$ in $G_R$.
\item Every period $c$ of rank greater than $R$ satisfies
\[
 |\can(e)|<\eta|c|
\]
for every $e\in E$.
\end{enumerate}
\end{lemma}
\begin{proof}
For each $e\in E$, the word $\can(e)^n$ is trivial in $G$, so it is a finite product of conjugates of relators from the sets $\mathrm{Rel}_i$ and their inverses. Only finitely many ranks occur in such an expression. Since $E$ is finite, there is a common stage $R_0\geq1$ at which all these relations hold.

We next explain why later periods become long. By the definitions in \cite[Sections~4.2--4.3]{ART}, rank-two periods have length at least two. A rank-three period must cyclically contain a $\tau$-th power of a nonempty word, otherwise it would already have rank two. Its length is therefore at least $\tau$. For $i\geq4$, a rank-$i$ period must cyclically contain $b^\tau$ for some period $b$ of rank $i-1$, or the preceding rank would already apply.
Induction therefore gives
\begin{equation}\label{eq:period-lengths}
 |a|\geq2,\qquad |c|\geq\tau^{i-2}
\end{equation}
for $a^n\in\mathrm{Rel}_2$ and $c^n\in\mathrm{Rel}_i$ with $i\geq3$.

Choose $M\geq0$ bounding the lengths of all words $\can(e)$ with $e\in E$. Now choose $R\geq\max\{R_0,2\}$ so large that $M<\eta\tau^{R-1}$.
For every $c^n\in\mathrm{Rel}_i$ with $i>R$, we then have
\[
 |\can(e)|\leq M
 <\eta\tau^{R-1}
 \leq\eta|c|
\]
for every $e\in E$. The torsion relations still hold because $G_R$ is a
quotient of $G_{R_0}$.
\end{proof}

We will use $\eta=1/4$ in the preceding lemma. The next lemma describes how periodic blocks survive canonical multiplication, with bounded endpoint corrections depending only on the adjacent coefficients and signs.

\begin{lemma}[Controlled multiplication]
\label{lem:uniform-buffering}
Fix $a^n\in\mathrm{Rel}_\rho$ with $\rho\geq2$, and put $s=\rho-1$ and $g=\can_s(a)$. Then $g$ has infinite order in $G_s$, and there is an
integer $K$ with $1-\tau\leq K\leq2\tau$ for which the following holds.

Let $\mathcal S\subseteq\mathrm{Can}_s$ be finite and closed under inverses, and suppose that every $c\in\mathcal S$ satisfies
\[
 c\notin\langle g\rangle\text{ in }G_s,\qquad
 \ell_\rho(c)<1,\qquad \ell_i(c)<\tfrac14
\]
for every $i>\rho$. Consider the canonical product
\[
 U=\can_s\bigl(c_0g^{e_1}c_1\cdots g^{e_t}c_t\bigr),
\]
where $t\geq1$, the $e_j$ are nonzero integers, $c_1,\ldots,c_{t-1}\in\mathcal S$, and $c_0,c_t\in\mathcal S\cup\{1\}$.
Put $\varepsilon_j=\operatorname{sign}(e_j)$ and $p_j=|e_j|-K$.
The shift $K$ accounts for the difference between the absolute exponent and the length of the normalized periodic part.

There are real-valued correction functions $\alpha,\beta$ giving the left and right endpoint corrections $\alpha_j,\beta_j$ for each block.
At interior ends, these are
\[
 \begin{aligned}
 \alpha_j&=\alpha(\varepsilon_{j-1},c_{j-1},\varepsilon_j),\\
 \beta_j&=\beta(\varepsilon_j,c_j,\varepsilon_{j+1}),
 \end{aligned}
\]
for $2\leq j\leq t$ and $1\leq j<t$, respectively. At the exterior ends, $\alpha_1$ depends only on $(c_0,\varepsilon_1)$, and $\beta_t$ only on $(\varepsilon_t,c_t)$. These corrections are chosen uniformly for all such products and satisfy
\begin{equation}\label{eq:correction-range}
 -3\tau-3<\alpha_j,\beta_j<3\tau+4
\end{equation}
for every $1\leq j\leq t$.

Suppose that
\[
 p_j\geq\tau,\qquad
 p_j+\alpha_j+\beta_j>\tau+1
\]
for every $1\leq j\leq t$. Then the following hold.
\begin{enumerate}
\item \emph{Surviving periodic blocks.}
The word $U$ contains distinct maximal $a$-periodic occurrences $v_1,\ldots,v_t$, in this order, with signs $\varepsilon_j$ and measures
\[
 \Lambda_a(v_j)=p_j+\alpha_j+\beta_j.
\]

\item \emph{Dependence on adjacent data.}
Each endpoint phase depends only on the same data as its correction. For $1\leq j<t$, the gap word between $v_j$ and $v_{j+1}$ depends only on $(\varepsilon_j,c_j,\varepsilon_{j+1})$. These data are independent of the exponent magnitudes and of all other junctions. A gap is taken to be empty when its adjacent occurrences overlap.

\item \emph{Other periodic occurrences.}
Let $J$ be a gap between consecutive occurrences, the prefix before $v_1$, or the suffix after $v_t$. Then
\[
 \ell_i(J)<4\tau+\frac{17}{4}
\]
for every $i>\rho$. Every maximal rank-$\rho$ periodic occurrence other than the $v_j$ has measure less than $8\tau+5$.
\end{enumerate}
\end{lemma}

The proof is given in Subsection~\ref{subsec:local-proof}. In our applications, the nonidentity coefficients have finite order in $G_s$, while $g$ has infinite order. This verifies the condition $c\notin\langle g\rangle$.

\subsection{Guard lengths and later periods}

Keep the notation of Lemma~\ref{lem:uniform-buffering} and fix $b\in\mathcal S$. The next construction places a central positive power of $g$ between two \emph{guards} $d_L,d_R$, separated from it by $b$ on each side. The guards must survive changes at the exterior junctions and remain shorter than the central marker, even when two adjacent guards merge.

The fixed inner junctions have data $(1,b,1)$. Put $\alpha=\alpha(1,b,1)$ and $\beta=\beta(1,b,1)$.
By the measure of a guard we mean the measure of its distinguished maximal $a$-periodic occurrence in the canonical product.

To leave more than $\tau+2$ periods after the variable endpoint correction, choose integers $p_L,p_R$ such that
\begin{equation}\label{eq:guard-calibration}
 4\tau+5<p_L+\beta\leq4\tau+6,\qquad
 4\tau+5<p_R+\alpha\leq4\tau+6,
\end{equation}
and put $d_L=g^{p_L+K}$ and $d_R=g^{p_R+K}$.
Such integers exist because each interval has length one, and \eqref{eq:correction-range} gives $p_L,p_R\geq\tau+2$.
The variable endpoint correction lies in $(-3\tau-3,3\tau+4)$, so each guard has measure in $(\tau+2,7\tau+10)$ in every permitted context.

Two guards meeting as $d_Rd_L$ between fixed $b$-neighbors merge into the power $g^{p_L+p_R+2K}$. Their parameter in the multiplication lemma is therefore $p_L+p_R+K\geq\tau+5$.
Using $1-\tau\leq K\leq2\tau$, the merged measure satisfies
\begin{equation}\label{eq:merged-guard}
 7\tau+11
 <
 (p_R+\alpha)+(p_L+\beta)+K
 \leq10\tau+12.
\end{equation}
Inverting the configurations gives the same bounds for negative powers, with $b^{-1}$ as the fixed coefficient.

Lemma~\ref{lem:finite-stage-data} gives a stage $R$ suited to our finite coefficient set. The next lemma ensures that we can choose a defining period after that stage.

\begin{lemma}[Periods beyond a prescribed stage]
\label{lem:unbounded-ranks}
For every integer $R\geq1$, there exist an integer $\rho>R$ and a word $a$ such that
\[
 a^n\in\mathrm{Rel}_\rho.
\]
\end{lemma}

\begin{proof}
The word $uv$ in two distinct basis letters is a rank-two period. We show that every rank $\rho\geq2$ with $\mathrm{Rel}_\rho\neq\varnothing$ is followed by another such rank.

Fix $a^n\in\mathrm{Rel}_\rho$ with $\rho\geq2$, and put $s=\rho-1$ and $g=\can_s(a)$. Choose a basis letter $b$. It has nontrivial finite order
in $G_s$, whereas $g$ has infinite order, so $b^{\pm1}\notin\langle g\rangle$. Together with \eqref{eq:period-lengths}, this verifies
the hypotheses of Lemma~\ref{lem:uniform-buffering} for $\mathcal S=\{b,b^{-1}\}$.

Choose the guards $d_L,d_R$ as in \eqref{eq:guard-calibration}, and for every $k\geq0$ put
\[
 U_k=\can_s\bigl(d_Lb(d_Rb)^k d_R\bigr).
\]
By Lemma~\ref{lem:uniform-buffering} and the guard bounds, $U_k$ contains $k+2$ distinguished periodic blocks of measure in $(\tau+2,7\tau+10)$, with consecutive blocks overlapping by less than one period.
Hence $|U_k|\to\infty$. The same estimates give
\[
 \ell_\rho(U_k)<8\tau+5<T,
\]
so $U_k\in\mathrm{Can}_\rho$ by IH 4.

It follows that $d_Rb$ has infinite order in $G_\rho$: otherwise, the products defining $U_k$ would take only finitely many values, contradicting the unbounded lengths of their canonical representatives.
Its $n$th power is nontrivial in $G_\rho$ but trivial in $G$. Since $G$ is obtained by imposing all the ART relators \cite[Corollary~8.7]{ART}, some $\mathrm{Rel}_i$ with $i>\rho$ must be nonempty. Iteration proves the claim.
\end{proof}

\subsection{Construction of the one word}
\label{subsec:one-word-proof}

We construct $w$ in Lemma \ref{lem:one-word-avoidance} so that every factor $w^{\pm1}$ leaves a recognizable marker in the canonical product.
The markers determine the signs, and the words between successive markers determine the intervening coefficients. We then show that a long later-rank power would force a forbidden repetition in this sign--coefficient record.

\begin{proof}[Proof of Lemma~\ref{lem:one-word-avoidance}]
\emph{Choose the guards and marker.}
Enlarge $E$ to be inverse-closed and contain some $b\ne1$. Apply Lemma \ref{lem:finite-stage-data} with $\eta=1/4$, and choose $a^n\in\mathrm{Rel}_\rho$ beyond the resulting stage by Lemma \ref{lem:unbounded-ranks}. Put $s=\rho-1$ and $g=\can_s(a)$.
The nonidentity chosen lifts of $E$ then satisfy Lemma \ref{lem:uniform-buffering}.

Use the guards $d_L,d_R$ calibrated for $b$ in \eqref{eq:guard-calibration}. With $\alpha=\alpha(1,b,1)$ and $\beta=\beta(1,b,1)$, choose an integer $Q$ such that
\[
 10\tau+13<Q+\alpha+\beta\leq10\tau+14.
\]
Then $Q\geq4\tau+6$. Set
\begin{equation}\label{eq:guarded-word}
 w=\underbrace{d_L}_{\text{left guard}}\,b\,
   \underbrace{g^{Q+K}}_{\text{marker}}\,b\,
   \underbrace{d_R}_{\text{right guard}}.
\end{equation}
Let $U$ be the rank-$s$ canonical form of a product $e_0w^{\varepsilon_1}e_1\cdots w^{\varepsilon_t}e_t$ as in the lemma.
At an identity interior coefficient the signs agree, so the adjacent guards merge as $d_Rd_L$ or its inverse. All remaining interior coefficients are nonidentity. Controlled multiplication and \eqref{eq:merged-guard} give
\begin{center}
\begin{tabular}{lcc}
\hline
Piece & Periodic measure & $\tau=15$\\
\hline
Guard, merged or unmerged
  & $(\tau+2,10\tau+12]$
  & $(17,162]$\\
Central marker
  & $(10\tau+13,10\tau+14]$
  & $(163,164]$\\
Other maximal rank-$\rho$ occurrence
  & $<8\tau+5$
  & $<125$\\
\hline
\end{tabular}
\end{center}
Thus the central markers are exactly the maximal rank-$\rho$ occurrences of measure $>10\tau+13$, so they can be recognized from the canonical word.

\emph{Read the signs and coefficients.}
For each sign $\varepsilon\in\{\pm1\}$, both junctions adjacent to the central block have data $(\varepsilon,b^\varepsilon,\varepsilon)$, so
Lemma \ref{lem:uniform-buffering} makes its endpoint corrections and phases, and hence the marker word $m_\varepsilon$, independent of the exterior coefficients, with $m_-=m_+^{-1}$ by inversion. They determine their signs, since a nonempty reduced word differs from its inverse. The markers are disjoint: each intervening guard has measure $>\tau+2$ and overlaps either neighboring block by less than one period.
Let $D(\varepsilon,e,\delta)$ be the connecting word for the indicated signs and coefficient. Locality gives the literal form
\begin{equation}\label{eq:marked-form}
 U=Lm_{\varepsilon_1}D(\varepsilon_1,e_1,\varepsilon_2)m_{\varepsilon_2}\cdots D(\varepsilon_{t-1},e_{t-1},\varepsilon_t)m_{\varepsilon_t}R.
\end{equation}
Write $\can_s(w^\varepsilon)=P_\varepsilon m_\varepsilon S_\varepsilon$.
The canonical form of $w^\varepsilon e w^\delta$ retains the prefix $P_\varepsilon m_\varepsilon$ and suffix $m_\delta S_\delta$, even when guards merge. Cancelling these ends gives
\begin{equation}\label{eq:one-word-gap}
 D(\varepsilon,e,\delta)=S_\varepsilon eP_\delta
 \quad\text{in }G_s.
\end{equation}
Thus each connecting word, together with its two signs, determines $e$ in $G_s$ and hence in $G$.

\emph{A later power would repeat the interior record.}
Let $v^n\in\mathrm{Rel}_i$ with $i>\rho$, and put $p=|v|$. A $v$-periodic interval containing no complete marker is covered by at most three junction gaps (including an exterior prefix or suffix), two guards, and two marker fragments. Indeed, meeting three markers would force the interval to contain the middle one completely; between consecutive markers there are at most two guards and three junction gaps, while before the first or after the last marker there are at most one guard and two gaps. Each gap contributes less than $(4\tau+17/4)p$ by Lemma \ref{lem:uniform-buffering}, and each remaining piece contributes less than $2p$ by IH 5. Its length is therefore less than
\begin{equation}\label{eq:marker-free-bound}
 \left[3\left(4\tau+\frac{17}{4}\right)+4\cdot2\right]p
 =\left(12\tau+\frac{83}{4}\right)p.
\end{equation}

Suppose that a $v$-periodic interval $u$ has length at least $(12\tau+21)p$. By \eqref{eq:period-lengths}, $p\geq\tau>8$, so the margin $p/4$ over \eqref{eq:marker-free-bound} exceeds two letters.
Deleting the endpoint letters of $u$ therefore leaves a complete marker, with both stopping letters inside $u$.

All markers have the same length, which is $<2p$ by IH 5 applied to the marker just found. Take the start of $u$ as coordinate zero. By periodicity, we can translate this marker by an integer multiple of $p$
to a start $x\in[2p,3p)$. Its copies at $x+jp$, for $0\leq j\leq\tau$, lie with both stopping letters inside $u$: the last ends before
$(\tau+5)p$. The stopping letters certify maximality, so the table
identifies all these copies as central markers.

For $0\leq j<\tau$, list the markers starting in
\[
 I_j=[x+jp,x+(j+1)p)
\]
in their order of occurrence, pairing each sign with the coefficient recovered from the following connecting word. The common length bound keeps these markers and their stopping letters inside $u$. Translation by $p$ preserves the markers and connecting words, so \eqref{eq:one-word-gap} gives identical lists of sign--coefficient pairs. Each list is nonempty because of the marker at $x+jp$; the marker at $x+\tau p$ supplies the last successor, so every recorded coefficient is interior.
These lists give $\tau$ consecutive copies of the same nonempty block in the interior record, contradicting its $\tau$-power-freeness.

We have proved
\[
 U\in\mathrm{Can}_s,\qquad
 \ell_\rho(U)\leq10\tau+14<T,\qquad
 \ell_i(U)<12\tau+21<T
\]
for every $i>\rho$. Successive applications of IH~4 leave $U$ unchanged at every later rank. Since $U$ contains a central marker, it is nonempty and represents a nonidentity element of $G$.
\end{proof}

\subsection{Proof of the controlled multiplication lemma}
\label{subsec:local-proof}

We first normalize the powers of $g$. At a junction, the condition $c\notin\langle g\rangle$ prevents cancellation of a full period between the periodic parts and prevents them from merging. The remaining estimates control their endpoints and the intervening gap.

\begin{lemma}
\label{lem:normalized-powers}
Fix $a^n\in\mathrm{Rel}_\rho$ with $\rho\geq2$, put $s=\rho-1$, and set $g=\can_s(a)$. There are words $W,Z$ and an integer $K$ such that
\begin{equation}\label{eq:marker-power}
 \can_s(g^{p+K})=Wa^pZ,\qquad
 1-\tau\leq K\leq2\tau
\end{equation}
for every integer $p\geq\tau$. Here $a$ is neither a suffix of $W$ nor a prefix of $Z$, and
\begin{equation}\label{eq:tail-bounds}
 \begin{aligned}
 \ell_a(W),\ell_a(Z)&<\tau+1,\\
 \ell_\rho(W),\ell_\rho(Z)&<2\tau+1,\\
 \ell_i(W),\ell_i(Z)&<\tau+2
 \end{aligned}
\end{equation}
for every $i>\rho$. Moreover, $g$ has infinite order in $G_s$, and
\begin{equation}\label{eq:marker-conjugacies}
 g=WaW^{-1}=Z^{-1}aZ,\qquad ZW=a^K
 \quad\text{in }G_s.
\end{equation}
\end{lemma}

\begin{proof}
By the definition of a rank-$\rho$ period, $a$ is cyclically canonical at rank $s$. Lemma~4.3 and Corollary~3.7 of \cite{ART} give fixed words $U,V$ such that
\[
 \can_s(g^M)=Ua^{M-2\tau}V
\]
for $M\geq4\tau$. The left end consists of a triangle side followed by a
terminal piece of $a^\tau$; the right end has the reversed description. Strip maximal whole powers, writing $U=Wa^{k_1}$ and $V=a^{k_2}Z$.
The remaining $W$ is a triangle-side prefix, possibly followed by less than one $a$-period; $Z$ has the symmetric form. The triangle side contributes less than $\tau$ in each later-rank measure. By IH 5, the remaining periodic piece contributes less than one for $a$, less than $\tau+1$ at rank $\rho$, and less than two at higher ranks. Subadditivity gives \eqref{eq:tail-bounds}.

Set $K=2\tau-k_1-k_2\leq2\tau$.
Then $\can_s(g^M)=Wa^{M-K}Z$ for all sufficiently large $M$. Comparing consecutive large powers gives the conjugacy identities in \eqref{eq:marker-conjugacies}, and substituting them into the power formula gives $ZW=a^K$. The unbounded canonical lengths show that $g$ has infinite order in $G_s$. IH 11 makes $Wa^pZ$ canonical for every $p\geq\tau$, and \eqref{eq:marker-conjugacies} identifies its value as $g^{p+K}$. If $K\leq-\tau$, taking $p=-K$ would give a nonempty canonical representative of the identity. Thus, $K\geq1-\tau$.
\end{proof}

For negative powers, use the inverted coordinates $(a^{-1},g^{-1},Z^{-1},W^{-1})$ with the same $K$. We call the displayed $a^p$ the \emph{periodic part}. It can extend into either end word by less than one period. Endpoint corrections are measured relative to this part, with positive values recording extensions and negative values recording losses.

\begin{proof}[Proof of Lemma~\ref{lem:uniform-buffering}]
We first analyze one junction between long powers and assemble these calculations into a canonical product. We then shorten the periodic cores, checking both the group identity and canonicality.

\emph{One junction.}
Fix $(\varepsilon,c,\delta)$. Using inverted coordinates if necessary, assume $\varepsilon=1$.
Put $q=4\tau+4$ and, for sufficiently large $P,P'$, write
\[
 A=\can_s(g^{P+K})=La^qZ,\qquad
 B=\can_s(g^{\delta(P'+K)})=H_\delta a^{\delta q}R,
\]
where $H_+=W$ and $H_-=Z^{-1}$. The loss bound below leaves more than
$q-(3\tau+3)=\tau+1$ reserved periods on each side, as needed for locality.
Apply canonical multiplication twice:
\begin{equation}\label{eq:two-local-triangles}
 U:=\can_s(Ac)=A'Dc',\qquad
 J:=\can_s(UB)=U'EB'.
\end{equation}
Here $A',U'$ are retained prefixes of $A,U$, respectively, $c',B'$ are retained suffixes of $c,B$, and $D,E$ are the inserted triangle sides.
Every triangle side has measure less than $\tau$ at each rank greater than $s$.

In the first multiplication, free cancellation removes less than one period from the periodic part, since $\ell_a(c)<1$, and the triangle removes less than $\tau$ more.
The resulting right endpoint correction lies in $(-\tau-1,\tau+2)$.
For the upper bound, a retained prefix of $Z$ contributes less than one period, since $Z$ does not start with $a$; $D$ and $c'$ contribute less than $\tau$ and one, respectively. The tail $Z_0$ after this maximal occurrence is covered by subwords of $Z,D,c$, so
\begin{equation}\label{eq:first-local-tail}
 \ell_a(Z_0)<(\tau+1)+\tau+1=2\tau+2.
\end{equation}

Suppose that a full period cancels directly between the surviving left periodic occurrence and the right periodic part. If $\delta=1$, this would identify cyclic shifts of $a$ and $a^{-1}$, which is impossible.
If $\delta=-1$, primitivity of $a$ aligns whole-period boundaries in the matched portion. Choose corresponding cuts, writing $U=Wa^hC$ and cutting the right factor after $Z^{-1}a^{-k}$. The intervening word cancels freely, so $\red(CZ^{-1}a^{-k})=1$.
Consequently, \eqref{eq:marker-conjugacies} gives
\[
 g^{P+K}c=Wa^{h+k}Z=g^{h+k+K}
 \quad\text{in }G_s,
\]
contrary to $c\notin\langle g\rangle$. Thus, the periodic parts cancel directly by less than one period.

The right part loses less than $2\tau+2$ periods against $Z_0$, less than one against the left periodic part, and less than $\tau$ to the second triangle. The left part first lost less than $\tau+1$; it then loses less than $\tau+1$ against $H_\delta$, less than one against the opposite periodic part, and less than $\tau$ to the second triangle.
Both total losses are therefore less than $3\tau+3$.

The surviving maximal occurrences cannot merge. Opposite signs have incompatible periodic continuations. For equal signs, a merger would give $J=Wa^uZ$ for some integer $u$, since the outer end words remain unchanged. Then
\[
 g^{P+K}c\,g^{P'+K}=g^{u+K}
 \quad\text{in }G_s,
\]
again forcing $c\in\langle g\rangle$.

For the left occurrence to extend beyond its first-stage endpoint, its old stopping letter must be removed.
Since $U'$ is a prefix of $U$, it then contains none of $Z_0$. Beyond the first correction, the extension can use only $E$, the remaining part of $H_\delta$, and less than one period of overlap with the right occurrence. Its upper bound is therefore
\[
 (\tau+2)+\tau+(\tau+1)+1=3\tau+4.
\]
Inverting the product gives the corresponding bound at the other end. This proves \eqref{eq:correction-range}.

\emph{Gap estimates and locality.}
Any nonempty gap is covered, in order, by surviving subwords of
\[
 Z,\quad D,\quad c,\quad E,\quad H_\delta.
\]
Indeed, \eqref{eq:two-local-triangles} retains prefixes and suffixes, and each remaining periodic part belongs to its maximal occurrence.
Thus
\begin{equation}\label{eq:local-gap-covers}
 \begin{aligned}
 \ell_\rho(\text{gap})
 &<2(2\tau+1)+1+2\tau=6\tau+3,\\
 \ell_i(\text{gap})
 &<2(\tau+2)+\tfrac14+2\tau
   =4\tau+\tfrac{17}{4}
 \end{aligned}
\end{equation}
for every $i>\rho$.
Overlapping occurrences have empty gap. At exterior ends, a single canonical multiplication gives the same correction and gap bounds using fewer pieces; an identity exterior coefficient is handled directly by
the normalized form.

Apply Corollary 3.9 of \cite{ART}, and its inverse, successively to the two multiplications. When an exponent is changed, the replacement power is canonical by Lemma \ref{lem:normalized-powers}, and the first locality application supplies the canonical intermediate input for the second. Thus, the endpoint corrections, phases and gap words depend only on $(\varepsilon,c,\delta)$. These define the functions $\alpha,\beta$ in the statement, with the analogous one-sided data at exterior ends.

\emph{Assembly for long powers.}
Keep $\varepsilon_j,p_j$ from the statement.
Choose auxiliary integers $P_j\geq p_j$ sufficiently large that the two reserved blocks of $q$ periods in each factor are disjoint, and consider
\[
 c_0g^{\varepsilon_1(P_1+K)}c_1\cdots
 g^{\varepsilon_t(P_t+K)}c_t.
\]
Apply the local calculation successively from left to right. At each step the preceding partial product is canonical, and a surviving block of $\tau$ periods separates the new junction from the earlier ones, so Corollary 3.9 applies. The resulting maximal occurrences have measures
$P_j+\alpha_j+\beta_j$.

From the left end of the $j$-th occurrence, select a core consisting of
\[
 N_j=\lfloor P_j+\alpha_j+\beta_j\rfloor-1
\]
complete periods. The subtraction of one leaves a right tail of measure
in $[1,2)$, making the cores disjoint because distinct maximal $a$-periodic occurrences overlap by less than one period.
The canonical word has the literal form
\begin{equation}\label{eq:long-core-form}
 A_0\widehat a_1^{N_1}A_1\cdots
 \widehat a_t^{N_t}A_t,
\end{equation}
where $\widehat a_j$ is the cyclic shift of $a^{\varepsilon_j}$ read at its left end. The fixed endpoint phases and tail lengths, together with
locality, make all $\widehat a_j$ and $A_j$ independent of the $P_j$.

\emph{Shortening the cores.}
First extend the group identity between \eqref{eq:long-core-form} and the auxiliary product. Fix all parameters except $P_j$.
Since $N_j-P_j$ is constant, comparison at two consecutive values gives
\[
 X\widehat a_jX^{-1}=Yg^{\varepsilon_j}Y^{-1}
 \quad\text{in }G_s,
\]
where $X,Y$ are the products of the factors preceding the corresponding powers in the two expressions. For $h\in\mathbb Z$, replacing $P_j$ by $P_j+h$ left-multiplies the two expressions by $X\widehat a_j^hX^{-1}$ and $Yg^{\varepsilon_j h}Y^{-1}$, respectively. These are equal by the displayed identity. Thus equality at one value extends to every integer
value of $P_j$. Proceeding one coordinate at a time extends the group identity to every integer tuple.

To check canonicality, replace each $P_j$ by $p_j$. The new core exponent is
\[
 \lfloor p_j+\alpha_j+\beta_j\rfloor-1\geq\tau.
\]
Each $\widehat a_j$ is a rank-$\rho$ period, so IH 11 permits these replacements successively, preserving canonicality and the intervening words. The extended group identity identifies the result as the canonical form of the required product.

The fixed right tails keep every left stopping letter in an intervening word $A_j$. A right stopping letter lies either in such a word or in
the first period of the next core, including when occurrences overlap or abut. That first period survives. Since only whole periods are removed, the endpoint phases are unchanged. Thus maximality, overlaps and gap words are preserved, and the resulting occurrences $v_j$ have measures $p_j+\alpha_j+\beta_j$.

Finally, a further maximal rank-$\rho$ occurrence cannot contain an entire $v_j$: compatible phases would contradict maximality, and incompatible phases would contradict IH 5 because $\Lambda_a(v_j)>\tau+1$. It is therefore covered by at most one gap or exterior region and two bordering block fragments. Its measure is less than
\[
 (6\tau+3)+2(\tau+1)=8\tau+5,
\]
as required.
\end{proof}

\Needspace{8\baselineskip}

\bibliographystyle{amsplain}
\bibliography{burnside_selflessness}
\end{document}